\documentclass{elsart}

\usepackage{graphicx}
\usepackage{amssymb}
\usepackage{epic,eepic}

\newcommand{\pe}[2]{\left\langle#1,#2\right\rangle}
\newcommand{\racion}[2]{\mbox{\small$\frac{{#1}}{{#2}}$}}

\begin{document}

\begin{frontmatter}

\title{Extensions of discrete classical orthogonal
polynomials beyond the orthogonality}
\author[RCS]{R. S. Costas-Santos\thanksref{RCST}}
\thanks[RCST]{The research of RSCS has been supported by
Direcci\'{o}n General de Investigaci\'{o}n del Ministerio de
Educaci\'{o}n y Ciencia of Spain under grant
MTM2006-13000-C03-02.}
\ead{rscosa@gmail.com}
\ead[url]{http://www.rscosa.com}
and
\author[JSL]{J. F. S\'{a}nchez-Lara\thanksref{JSLT}}
\ead{jlara@ual.es}
\thanks[JSLT]{The research of JFSL has been supported by
Secretar\'{\i}a de Estado de Universidades e Investigaci\'{o}n del
Ministerio de Educaci\'{o}n y Ciencia of Spain and by Junta de
Andaluc\'{\i}a, grant FQM229.}
\address[RCS]{Department of Mathematics, University of
California, South Hall, Room 6607 Santa Barbara, CA 93106 USA}
\address[JSL]{Universidad Polit\'{e}cnica de Madrid.
Escuela T\'{e}cnica Superior de Arquitectura.
Departamento de Matem\'{a}tica Aplicada.
Avda Juan de Herrera, 4.
28040 Madrid, Spain}

\begin{abstract}
It is well-known that the family of Hahn polynomials
$\{h_n^{\alpha,\beta}(x;N)\}_{n\ge 0}$ is orthogonal
with respect to a certain weight function up to degree
$N$.
In this paper we prove, by using the three-term
recurrence relation which this family satisfies,
that the Hahn polynomials can be characterized
by a $\Delta$-Sobolev orthogonality for every $n$
and present a factorization for Hahn polynomials
for a degree higher than $N$.

We also present analogous results for dual Hahn,
Krawtchouk, and Racah polynomials and give the limit
relations among them for all $n\in \mathbb{N}_0$.
Furthermore, in order to get these results for the
Krawtchouk polynomials we will obtain a more general
property of orthogonality for Meixner polynomials.
\end{abstract}

\begin{keyword}
Classical orthogonal polynomials \sep inner product
involving difference operators \sep non–standard
orthogonality.

\emph{2000 MSC:} 33C45 \sep 42C05 \sep 34B24
\end{keyword}
\end{frontmatter}

\section{Introduction}

In the last decade some of the classical orthogonal
polynomials with non-classical parameters have
been provided with certain non-standard orthogonality.
For instance, K. H. Kwon and L. L. Littlejohn, in
\cite{kwli1}, established the orthogonality of the
generalized Laguerre polynomials $\{L_n^{(-k)}\}_{n
\ge 0}$, $k\ge 1$, with respect to the Sobolev inner
product:
$$
\pe f g=(f(0),f'(0),\dots,f^{(k-1)}(0))A\left(
\begin{array}{c} g(0)\\ g'(0)\\ \vdots \\ g^{(k-1)}(0)
\end{array}\right)+\int_0^{\infty}f^{(k)}(x)g^{(k)}(x)
e^{-x}dx,
$$
with $A$ being a symmetric $k\times k$ real matrix.
In \cite{kwli2}, the same authors showed that the Jacobi
polynomials $\{P_n^{(-1,-1)}\}_{n\ge 0}$, are orthogonal
with  respect to the inner product
$$
(f,g)_1=d_1f(1)g(1)+d_2 f(-1)g(-1)+ \int_{-1}^1
f'(x)g'(x)dx,
$$
where $d_1$ and $d_2$ are real numbers.

Later, in \cite{pepi}, T. E. P\'{e}rez and M. A. Pi\~{n}ar gave a
unified approach to  the orthogonality of the generalized
Laguerre polynomials $\{L_n^{(\alpha)}(x)\}_{n\ge 0}$, for
any real value of the parameter $\alpha$, by proving their
orthogonality with respect to a non-diagonal Sobolev  inner
product, whereas in \cite{pepi1} they have shown how to use
this orthogonality to obtain different properties of the
generalized Laguerre polynomials.

M. Alfaro, M.L. Rezola, T.E. P\'{e}rez and M.A. Pi\~{n}ar have
studied in \cite{alpepire} sequences of polynomials
which are orthogonal with respect to a Sobolev bilinear
form defined by
\begin{equation}
\label{bsq} {\cal B}_S^{(N)}=(f(c),f'(c),\dots,f^{(N-1)}
(c))A\left(\begin{array}{c}g(c)\\ g'(c) \\ \vdots \\
g^{(N-1)}(c)\end{array} \right)+\pe {{\bf u}} {f^{(N)}
g^{(N)}},
\end{equation}
where $\bf u$ is a quasi–definite linear functional, $c
\in \mathbb{R}$, $N$ is a positive integer, and $A$ is a
symmetric $N\times N$ real matrix such that each of its
principal submatrices is regular.

In particular, they deduced that Jacobi polynomials
$\{P^{(-N,\beta)}\}_{n\ge 0}$, where $N+\beta$ is not a
negative integer, are orthogonal with respect to
the bilinear form given in (\ref{bsq}), for $\bf u$ the
Jacobi functional corresponding to the weight function
$\rho^{(0,\beta+N)} (x)=(1+x)^{\beta+N}$ and $c=1$.

The remaining cases for the Jacobi polynomials, i.e. both
parameters, $\alpha$ and $\beta$, negative integers, were
considered by M. Alfaro, M. \'{A}lvarez de Morales and M.L.
Rezola in \cite{almore} where they proved that such
sequences satisfy a Sobolev orthogonality.

M. \'{A}lvarez de Morales, T. E. P\'{e}rez and M. A. Pi\~{n}ar  in
\cite{mopepi} have studied the sequence of the monic
Gegenbauer polynomials $\{C_n^{(-N+1/2)}\}_{n\ge 0}$,
where $N$ is a positive integer.
They have shown that this sequence is orthogonal with
respect to a Sobolev inner product of the form
$$
(f,g)_S^{2N}=(F(1)|F(-1))A(G(1)|G(-1))^T+
\int_{-1}^1f^{(2N)}(x)g^{(2N)}(x)(1-x^2)^Ndx,
$$
where
$$
(F(1)|F(-1))=(f(1),f'(1),\cdots,
f^{(N-1)}(1),f(-1),f'(-1),\cdots,f^{(N-1)}(-1)),
$$
$A=Q^{-1}D(Q^{-1})^T$, $Q$ is a non-singular matrix whose
entries are the consecutive derivatives of the Gegenbauer
polynomials evaluated at the points $1$ and $-1$, and $D$
is an arbitrary diagonal positive definite matrix.

M. \'{A}lvarez de Morales, T. E. P\'{e}rez, M. A. Pi\~{n}ar and A.
Ronveaux in \cite{mopepiro} have studied the sequence of
the monic Meixner polynomials $\{M_n^{(\gamma,\mu)}\}_{n
\ge 0}$, for $0<\mu<1$ and $\gamma\in \mathbb{R}$.
They have shown that this sequence is orthogonal with
respect to the inner product
$$
(f,g)_S^{(K,\gamma+K)}=\sum_{x=0}^{\infty }F(x)\Lambda^{(
K)}G(x)^T\rho^{(\gamma+K,\mu)}(x), \quad x\in [0,\infty),
$$
where $K$ is a non-negative integer, $F(x)=(f(x), \Delta
f(x), \cdots,\Delta^{K}f(x))$, $\Delta$ and $\nabla$  are,
respectively, the finite forward and backward difference
operators defined by
$$
\Delta f(x)=f(x+1)-f(x),\qquad \nabla f(x)=f(x)-f(x-1),
$$
$\rho^{(\gamma+K,\mu)}$ denotes the weight function
associated with the monic classical  Me\-ix\-ner
polynomials $\{M_n^{(\gamma+K,\mu)}\}_{n\ge 0}$, and
$\Lambda^{(K)}$ is a positive definite $(K+1)\times(K +1)$
matrix, with $K\ge \max\{0,-\gamma+1\}$.
Usually, when the inner product is defined by using the
difference operator instead the differential operator, the
orthogonality is said to be of $\Delta$-Sobolev type.

These examples suggest that the classical orthogonal
polynomials with non-classical parameters can be
provided with a Sobolev or a $\Delta$-Sobolev
orthogonality.
Furthermore, as was pointed out in \cite{mamo}, a
Sobolev-Askey tableau should be established.

In this paper we study discrete classical orthogonal
polynomials which exist only up to a certain degree.
This happens for the Krawtchouk, Hahn, dual Hahn and
Racah polynomials which satisfy a discrete orthogonality
with a finite number of masses.
These families exhibit this finite character in several
ways since there is a negative integer as a denominator
parameter in the hypergeometric representation.
Also in the three-term recurrence relation
$$
xp_n=p_{n+1}+\beta_np_n+\gamma_np_{n-1},\qquad n=0, 1, 2,
\dots
$$
there exists $M$ such that $\gamma_M=0$.
But there are other ways to characterize the discrete
classical orthogonal polynomials which, apparently, do
not say anything about whether the sequence of polynomials is
finite or infinite (see, for example, \cite{alal}).

We show that these polynomials can be considered for
degrees higher than $M$, in fact for all degrees,
and the main properties still hold except the orthogonality.
However, using difference properties of the polynomials,
we find a $\Delta$-Sobolev orthogonality of the form
$$
\pe f g_S:=\pe {{\bf u}_0}{fg}+\pe {{\bf
u}_1}{(\Delta^Mf)(\Delta^Mg)},
$$
with respect to which the polynomials are characterized,
where ${\bf u}_0$ and ${\bf u}_1$ are certain classical
functionals.

Also we obtain a factorization for these polynomials of
the form
$$
p_n=p_Mq_{n-M},\qquad n=M, M+1, M+2, \dots,
$$
where $p_M$ vanishes in the $M$ masses of the orthogonality
measure associated with the linear functional ${\bf u}_0$,
and $q_{n-M}$ is a classical orthogonal polynomial that is
at the same level in the Askey tableau as $p_n$.

The structure of the paper is as follows. The case of Hahn
polynomials is studied in Section 2 in detail.
In Sections 3, 4, and 5, we get analogous results for the
Racah, dual Hahn and Krawtchouk polynomials, which satisfy
a discrete orthogonality with a finite number of masses,
applying analogous reasoning that we have considered for
Hahn polynomials.
Finally in Section 6 we show that the known limit relations
involving the above families hold for any $n\in
\mathbb{N}_0$.
The Appendix is devoted to proving more general orthogonal
relations for Meixner polynomials which are used in Section
\ref{SecKrawtchouk}.

\section{Hahn polynomials}
The monic classical Hahn polynomials
$h_n^{\alpha,\beta}(x;N)$, $n=0, 1$, $\dots$, $N$, $N\in
\mathbb{N}_0$, can be defined by their explicit
representation in terms of the hypergeometric function
(see e.g. \cite[p.33] {kost}):
{\small \begin{equation}\label{HahnHipergeometric}
\hspace{-7mm} h_n^{\alpha,\beta}(x;N)\!=\!\frac{(-N,\alpha+
1)_n}{(\alpha\!+\!\beta\!+\!n\!+\!1)_n}{}_3F_2\!\!\left(
\left.\begin{array}{c}-n,\alpha+\beta+n+1,-x\\-N,\alpha+
1\end{array}\right|1 \right), \ n=0, 1, \dots,N,
\end{equation}}
where $(a)_n$ denotes the Pochhammer symbol
$$
(a)_0:=1,\qquad (a)_n:=a(a+1)\cdots (a+n-1),\qquad n=1, 2, 3,
\dots
$$
and $(a_1,a_2,\dots,a_j)_n:=(a_1)_n(a_2)_n\cdots(a_j)_n$.
These polynomials satisfy the following property of orthogonality:
\begin{equation}\label{HahnOrthogonality1}
\sum_{x=0}^{N}h_n^{\alpha,\beta}(x;N)h_m^{\alpha,\beta}(x;
N)\rho^{\alpha,\beta}(x;N)=0,\qquad 0\le m<n\leq N,
\end{equation}
where
$$
\rho^{\alpha,\beta}(x;N)=\frac{\Gamma(\beta+N+1-x)
\Gamma(\alpha+1+x)}{\Gamma(1+x)\Gamma(N+1-x)}.
$$
When $\alpha,\beta>-1$ or $\alpha,\beta<-N$ this is a
positive definite orthogonality. However, it is possible
to consider general complex parameters $\alpha$ and $\beta$
and (\ref{HahnOrthogonality1}) remains by using analytic
continuation.

Furthermore, Atakishiyev and Suslov \cite{atsu3} considered
Hahn polynomials for general complex parameters $\alpha,
\beta$ and $N$ and nowadays they are known as continuous
Hahn polynomials \cite{ask1}.
The monic ones are
$$
p_n(x;a,b,c,d)=D_n~ {}_3F_2\left(\left.\begin{array}{c}
-n,n+a+b+c+d-1,a+ix\\ a+c,a+d \end{array}\right|1\right),
$$
where
$$
D_n=\frac{i^n(a+c,a+d)_n}{(n+a+b+c+d-1)_n},
$$
and since the parameter $a$ causes only a translation, Hahn
and continuous Hahn polynomials are related in the
following way
\begin{eqnarray}
&p_n(x;a,b,c,d)=i^nh_n^{a+d-1,b+c-1}(-a-ix;-a-c),\nonumber
\\ &h_n^{\alpha,\beta}(x;N)=(-i)^np_n(ix;0,\beta+N+1,-N,
\alpha+1).\label{HahnToCHahn}
\end{eqnarray}
Continuous Hahn polynomials satisfy a non-Hermitian
orthogonality
\begin{equation}\label{ContHahnOrthogonality}
\int_Cp_n(x;a,b,c,d)x^mw(x;a,b,c,d)dx=0,\qquad 0\le m<n,
\end{equation}
where
\begin{equation} \label{weigthch}
w(x;a,b,c,d)=\Gamma(a+ix)\Gamma(b+ix)\Gamma(c-ix)
\Gamma(d-ix),
\end{equation}
and $C$ is a contour on $\mathbb{C}$ from $-\infty$ to
$\infty$ which separates the increasing poles
$$
(a+k)i, \ (b+k)i,\qquad k=0, 1, 2, \dots,
$$
from the decreasing ones
$$
(-c-k)i, \ (-d-k)i,\qquad k=0, 1, 2, \dots,
$$
which can be done when these two sets of poles are
disjoint, i.e.
$$
-a-c, \ -a-d,\ -b-c, -b-d \notin \mathbb{N}_{0}.
$$
The continuous Hahn polynomials also satisfy, among others,
a second order linear difference equation, a Rodrigues
formula, the TTRR which will be useful later
$$
xp_n(x)=p_{n+1}(x)+(B_n+a)ip_n(x)+
C_np_{n-1}(x),
$$
with
$$ \begin{array}{rl}
B_n=&\displaystyle \frac{n(n+b+c-1)(n+b+d-1)} {(2n+a+b+c+
d-2)(2n+a+b+c+d-1)} \\ & - \displaystyle  \frac{(n+a+b+c+
d-1)(n+a+c)(n+a+d)}{(2n+a+b+c+d-1)(2n+a+b+c+d)}
\end{array}$$
\begin{eqnarray}
C_n=&\frac{n(n+b+c-1)(n+b+d-1)(n+a+b+c+d-2)}{(2n+a+b+c+d-1)
(2n+a+b+c+d-2)}\nonumber \\ & \times \frac{(n+a+c-1)(n+a+d-
1)}{(2n+a+b+c+d-2)(2n+a+b+c+d-3)},
\label{coefRecCn}
\end{eqnarray}
and they have several generating functions.

Let us focus our attention on (\ref{HahnHipergeometric}), it
can be rewritten as
\begin{equation} \label{HahnHypergeometricreduced}
\hspace{-7mm} h_n^{\alpha,\beta}(x;N)\!=\! \frac{(\alpha+
1)_n}{(\alpha\!+\!\beta\!+\!n\!+\!1)_n}
\sum_{k=0}^n\frac{(-n,\alpha+\beta+n+1,-x)_k(-N+k)_{n-k}}
{(\alpha+1,1)_k},
\end{equation}
which is valid for every $n$ and, in this way, it
can be used to define Hahn polynomials for any $n\in
\mathbb{N}_0$.

These new polynomials satisfy the following result:

\begin{thm}
Let $N$ be a non-negative integer,
then the Hahn polynomials $h_n^{\alpha,\beta}(x;N)$ for
$n\ge 0$ satisfy the following properties:
\begin{itemize} \item[i)]
\begin{equation}\label{HahnRec1}\hspace{-7mm} \begin{array}
{l}h_{-1}^{\alpha,\beta}(x;N)=0,\quad  h_0^{\alpha,\beta}(x
;N)=1,\\ x  h_n^{\alpha,\beta} (x;N)= h_{n+1}^{\alpha,\beta}
(x;N)\!+\! \beta_n  h_n^{\alpha,\beta}(x;N)\!+\! \gamma_n
h_{n-1}^{\alpha,\beta}(x;N), \ n=0, 1, 2, \dots
\end{array} \end{equation}
where
\begin{eqnarray}
\beta_n=&\hspace{-12mm}\frac{(\alpha+1)N(\alpha+\beta)+n(2N-
\alpha+\beta)(\alpha+\beta+n+1)}{(\alpha+\beta+2n)(\alpha+
\beta+2n+2)}, \nonumber
\\
\gamma_n=&\frac{n(N\!+\!1-n)(\alpha\!+\!\beta+n)(\alpha+n)
(\beta+n)(\alpha+\beta+N\!+\!n\!+\!1)}{(\alpha+\beta+2n-1)
(\alpha+\beta+2n)^2(\alpha+\beta+2n+1)}.\label{coefRecGamma}
\end{eqnarray}
\item[ii)] For any integer $k$, $0\le k\le n$,
\begin{eqnarray}
ii.1) & \Delta^k  h_n^{\alpha,\beta}(x;N)=(n-k+1)_k
h_{n-k}^{\alpha+k,\beta+k}(x;N-k),\label{relii1}\\
ii.2) & \hspace{6mm} \nabla^k h_n^{\alpha,\beta}(x; N)=
(n-k+1)_k h_{n-k}^{\alpha+k,\beta+k}(x-k;N-k)\label{relii2}.
\end{eqnarray}
\item[iii)] Second order linear difference equation:
\begin{eqnarray}
& x(\beta+N+1-x)\nabla\triangle h_n^{\alpha,\beta}(x;N)+
\bigg((\alpha+1)N \\ & -(\alpha+\beta+2)x\bigg)\triangle
h_n^{\alpha,\beta}(x;N)+\lambda_n  h_n^{\alpha,\beta}(x;N)
=0, \label{HahnDifference}
\end{eqnarray}
with $\lambda_n=n(n+\alpha+\beta+1)$.
\item[vi)] Rodrigues formula:
\begin{equation}\label{HahnRodrigues}
h_n^{\alpha,\beta}(x;N)
\rho^{\alpha,\beta}(x;N)=\frac{(-1)^n}{(\alpha+\beta+n+
1)_n}\nabla^n(\rho^{\alpha+n,\beta+n}(x;N-n)),
\end{equation}
with
$$
\rho^{\alpha,\beta}(x;N)=\frac{\Gamma(\beta+N+1-x)
\Gamma(\alpha+1+x)}{\Gamma(1+x)\Gamma(N+1-x)}.
$$
\item[v)] Generating function:
\begin{equation}\label{HahnGenerating} \hspace{-7mm}
\sum_{n=0}^{\infty}\frac{(\alpha+\beta+n+1)_n}{(\beta+1,
\alpha+1,1)_n} h_n^{\alpha,\beta}(x;N)t^n \!\!=\!\!{}_1F_1
\!\!\left(\left.\begin{array}{c}-x\\ \alpha+1 \end{array}
\right|-t\right)\!\!{}_1F_1\!\!\left(\left. \begin{array}
{c} -N+x\\ \beta+1 \end{array}\right|t\right),
\end{equation}
valid for $x\in\mathbb{C}$ and
$-\alpha-1,-\beta-1\notin\mathbb{N}$.
\end{itemize}
\end{thm}

\begin{rem} \label{nota1} Apparently, the conditions
$-\alpha, -\beta,-\alpha-\beta\notin \mathbb{N}$ are
necessary in i)-iv), but this problem disappears by
using a suitable normalization, for instance, if there
is a polynomial dependence on the parameters.
\end{rem}

The proof is straightforward using the well-known
properties of continuous Hahn polynomials (see \cite{kost})
and the limit relation
$$
h_n^{\alpha,\beta}(x;N)=\lim_{\varepsilon\to 0}
(-i)^n~p_n(ix;0,\beta+N+\varepsilon+1,-N-\varepsilon,
\alpha+1),
$$
easily obtained from (\ref{HahnToCHahn}).
Note that for small $\varepsilon$ the weight function
associated with the continuous Hahn polynomials,
$w(x; 0, \alpha+N+ \varepsilon+1, -N- \varepsilon,
\beta+1)$, satisfies the poles separation condition.

Now we center our attention on the behavior of the zeros.
Figure 1 shows the standard configuration of the zeros
of Hahn polynomials for degree greater than $N+1$.
In fact, they have $N+1$ zeros on $0, 1,\dots,N$, and the
other $n-N-1$ zeros are located on unknown curves of the
complex plane.
In the special case $\alpha=\beta\in\mathbb{R}^+$ this
curve is the line $\{ti + N/2 : t\in\mathbb{R}\}$.
\begin{figure}[!hbt]
\center
\includegraphics[scale=0.35]{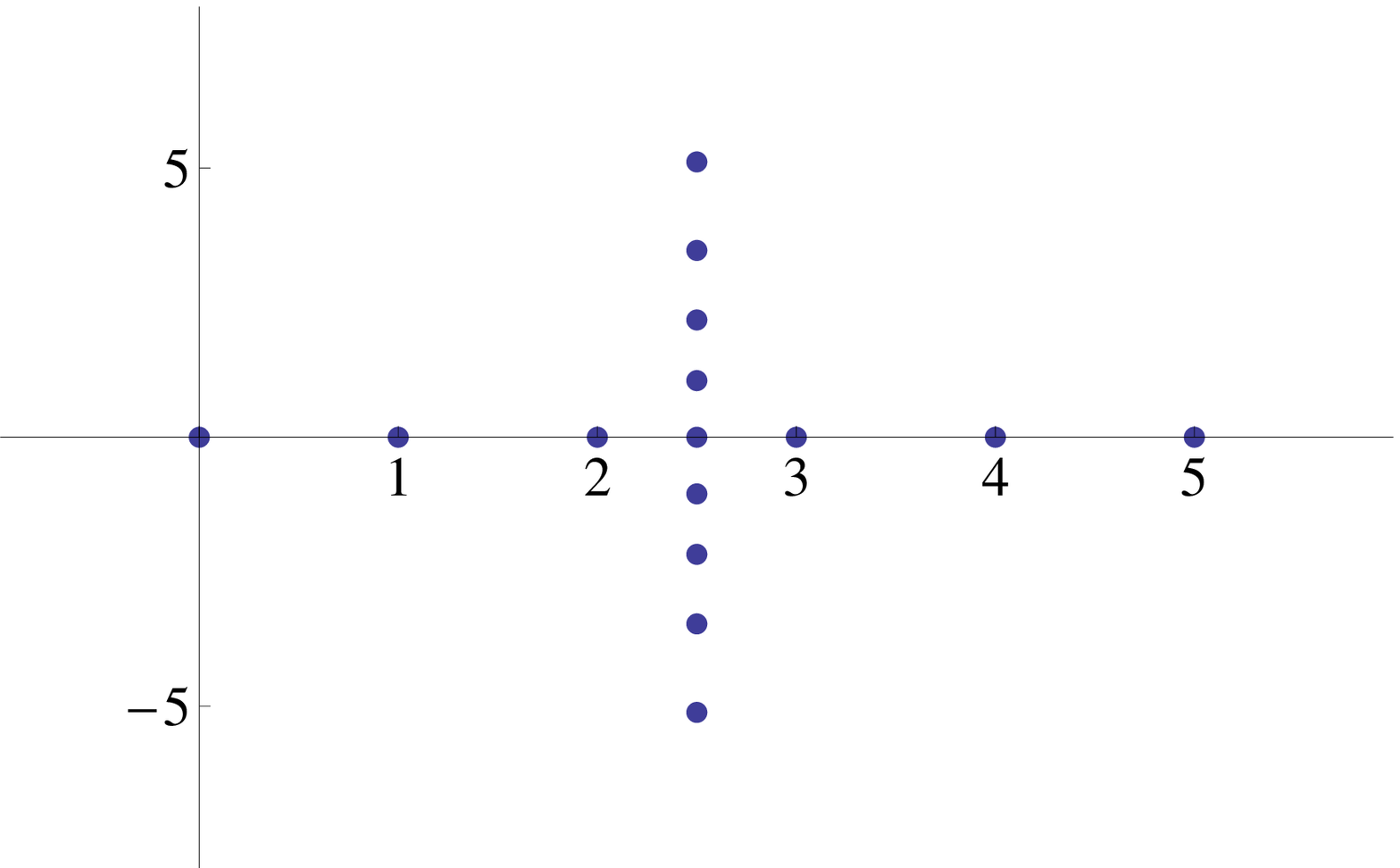}\hspace{1cm}
\includegraphics[scale=0.35]{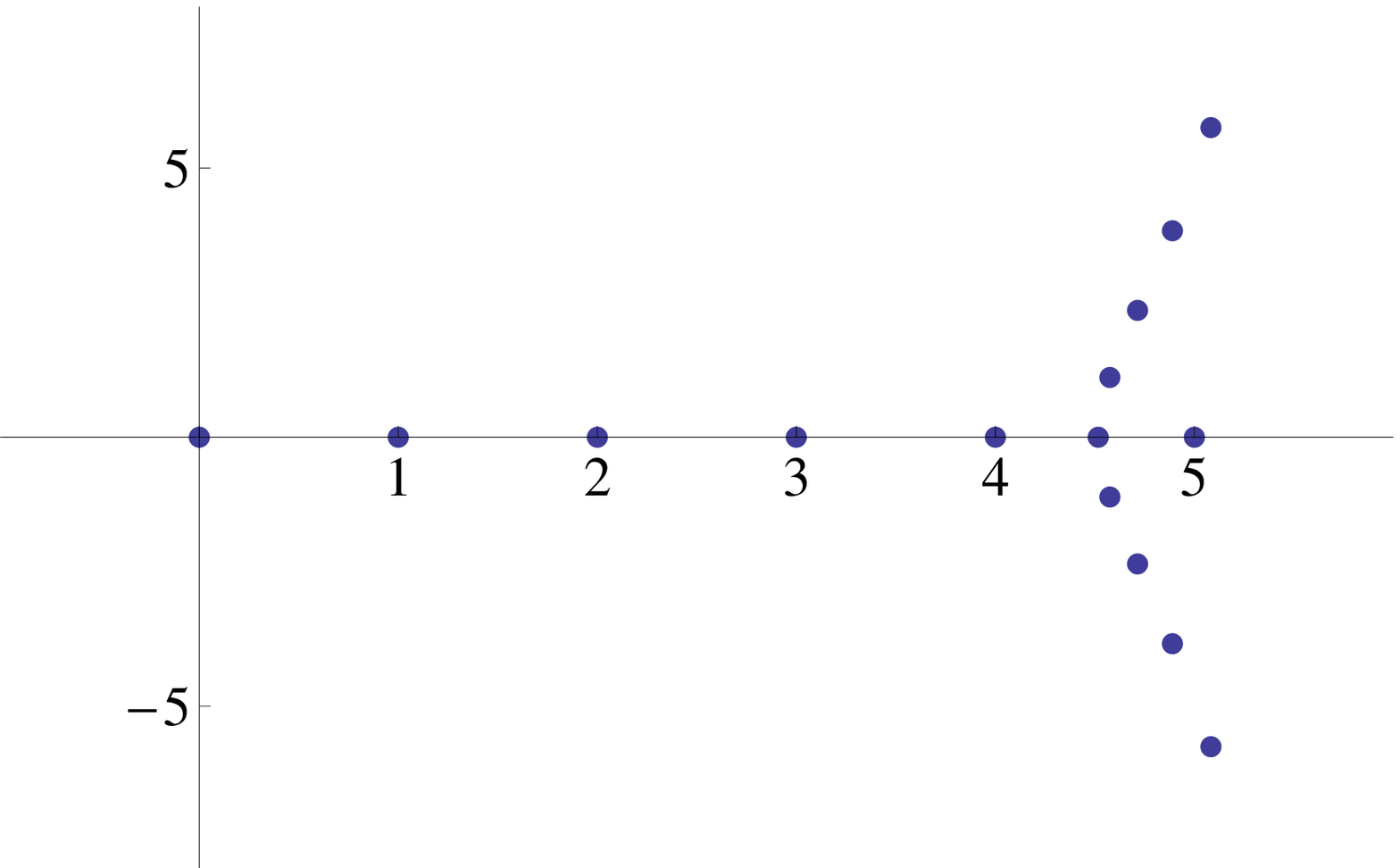} \caption{Zeros
of $h_{15}^{1,1}(x;5)$ (left) and $h_{15}^{1,15}(x;5)$
(right)} 
\end{figure}

The following result is straightforward taking into account
(\ref{HahnHypergeometricreduced}) and that
if $n\geq N+1$, then $(-N+k)_{n-k}=0$, for  $k=0,\dots,N$.

\begin{prop} \label{prop21} Let $N$ be a non-negative
integer.
For every $n\ge N+1$,
{\small $$
\begin{array}{rl}
h_n^{\alpha,\beta}(x;N)&=(x-N)_{N+1}(-i)^{n-N-1}p_{n-N-1}
(ix;N+1,\beta+N+1,1,\alpha+1)\\&\hspace{-20mm}=(x-N)_{N+1}
(-i)^{n-N-1}p_{n-N-1}\left(\left(x-\frac N2\right)i; 1+
\frac N2,\beta+1+\frac N2,1+\frac N2,\alpha+1+\frac N2
\right).
\end{array}
$$}
\end{prop}
\begin{rem} \label{nota2} Note that, as was expected
from (\ref{HahnOrthogonality1}), $h_{N+1}^{\alpha,\beta}
(x;N)=(x-N)_{N+1}$ which vanishes on $0, 1, \dots, N$ and
is independent of $\alpha$ and $\beta$ --- this fact is
non-trivial from the other ways to characterize these polynomials
(see e.g. \cite{alv3}).
On the other hand, in the case $\alpha=\beta$ the
continuous Hahn polynomial in the right-hand side is a
linear transformation of a real polynomial.
\end{rem}

Now we establish the main result.

\begin{thm} \label{theo22}
Let $N$ be a non-negative integer and $\alpha, \beta \in
\mathbb{C}$ such that
\begin{equation} \label{condortho}
-\alpha,-\beta \not \in \{1,2,\dots,N,N+2,\dots\},
\end{equation}
and
\begin{equation} \label{condortho-1}
-\alpha-\beta\not \in\{1, 2, \dots,
2N+1, 2N+3, \dots\}.
\end{equation}
Then the family of monic Hahn polynomials
$h_n^{\alpha,\beta}(x;N)$ for $n\ge 0$ is a MOPS with
respect to the following $\Delta$-Sobolev inner product:
$$
(f,g)_S=\sum_{x=0}^{N}f(x)g(x)\rho^{\alpha,\beta}(x;N)+
\int_{C}(\Delta^{N+1}f(z))(\Delta^{N+1}g(z))
\omega^{\alpha,\beta}(z;N)dz,
$$
where
$$
\begin{array}{rl}
\rho^{\alpha,\beta}(x;N)&=\displaystyle \frac{\Gamma(\beta+
N+1-x)\Gamma(\alpha+x+1)}{\Gamma(N+1-x)\Gamma(x+1)},\\
\omega^{\alpha,\beta}(z;N)&=\displaystyle \Gamma(-z)\Gamma(
\beta+N+1-z)\Gamma(1+z)\Gamma(\alpha+N+2+z),
\end{array}
$$
and $C$ is a complex contour  from $-\infty i$ to $\infty
i$ which separates the poles of the functions $\Gamma(-z)
\Gamma(\beta+N+1-z)$ and $\Gamma(1+z)\Gamma(\alpha+N+2+z)$.
Furthermore, this $\Delta$-Sobolev inner product
characterizes the polynomials $h_n^{\alpha,\beta}(x;N)$ for
all $n\in\mathbb{N}_0$.
\end{thm}

\begin{rem} Note that the conditions (\ref{condortho})
and (\ref{condortho-1}) are equivalent to the existence
of a unique $n$ such that $\gamma_n=0$, and therefore $n=N+1$.
Furthermore, (\ref{condortho}) is equivalent to $h_n^{
\alpha,\beta}(x;N)$ is uniquely determined by
(\ref{HahnOrthogonality1}) for $n\leq N$ together with the
poles separation condition of $w^{\alpha,\beta}(z;N)$ given
in the above theorem.
\end{rem}
\begin{pf} If $0\le m<n\le N+1$, since $\Delta^{N+1}x^m=0$,
$$
(h_n^{\alpha,\beta}(x;N),x^m)_S= \sum_{x=0}^{N}h_n^{
\alpha,\beta}(x;N)x^m\rho^{\alpha,\beta}(x;N)=0.
$$
If $n\ge N+1$ and $m\le N$, by Proposition \ref{prop21},
$$
(h_n^{\alpha,\beta}(x;N),x^m)_S=\sum_{x=0}^{N}h_n^{
\alpha,\beta}(x;N)x^m\rho^{\alpha,\beta}(x;N)=0,
$$
and if $N+1\le m<n$, then
$$
\begin{array}{rl}
(h_n^{\alpha,\beta}(x;N),x^m)_S=& \displaystyle \int_{C}
(\Delta^{N+1}h_n^{\alpha,\beta}(z;N))
(\Delta^{N+1}z^m)\omega^{\alpha,\beta}(z;N)dz\\[3mm]
= & \displaystyle \int_{C} p_{n-N-1}(zi;0,\beta+N+1,
1,\alpha+N+2)\\ & \times q_{m-N-1}(z)\omega^{\alpha,\beta}
(z;N)dz=0, \end{array}
$$
where $q_{m-N-1}$ is a polynomial of degree $m-N-1$.

Now we show that the inner product characterizes the
polynomials.
If $n\le N$ then, due to (\ref{condortho}), we get
$$
(h_n^{\alpha,\beta}(x;N),x^n)_S=\sum_{x=0}^{N}h_n^{\alpha,
\beta}(x;N) x^n\rho^{\alpha,\beta}(x;N)\ne 0,
$$
and if $n\ge N+1$ we obtain
$$
\begin{array}{rl}
(h_n^{\alpha,\beta}(x;N),x^n)_S= & \displaystyle \int_{C}
(\Delta^{N+1}h_n^{\alpha,\beta}(z;N))
(\Delta^{N+1}z^n)\omega^{\alpha,\beta}(z;N)dz\\[3mm]
= & \displaystyle \int_{C} p_{n-N-1}(zi;0,\beta+N+1,1,
\alpha+N+2)\\ & \times q_{n-N-1}(z)\omega^{\alpha,\beta}
(z;N)dz\ne 0, \end{array}
$$
since $q_{n-N-1}$ is a polynomial of degree $n-N-1$ and the
coefficient $C_k$, given in (\ref{coefRecCn}), is different
from zero for $a=0$, $b=\beta+N+1$, $c=1$, and $d=\alpha+N+
2$, for $k=1, 2, \dots$
\qed\end{pf}

\section{Racah polynomials}
We can apply an analogous process for the Racah polynomials
{\small
\begin{equation} \label{ex-Racah}
R_n(\lambda(x);\alpha,\beta, \gamma,\delta)=r_n \ {}_4F_3
\left.\left(\begin{array}{c}-n,n+\alpha+\beta+1,-x,x+\gamma+
\delta+1\\\alpha+1,\beta+\delta+1,\gamma+1\end{array}\right|
1 \right),
\end{equation}
}
in which
$$
r_n=\frac{(\alpha+1,\beta+\delta+1,\gamma+1)_n}{(n+\alpha+
\beta+1)_n},
$$
and $ \lambda(x)=x(x+\gamma+\delta+1)$, by using the Wilson
polynomials \cite[p. 28]{kost}
$$
W_n(x^2;a,b,c,d)=\ w_n \ {}_4F_3\left.\left(\!\!
\begin{array}{c}-n,n+a+b+c+d-1,a+ix,a-ix\\a+b,a+c,a+d
\end{array}\right|1\right),
$$
in which
$$
w_n=\frac{(-1)^n(a+b,a+c,a+d)_n}{(n+a+b+c+d-1)_n}.
$$
In fact,
$$
\begin{array}{rl}
&R_n(\lambda(x);\alpha,\beta,\gamma,\delta)\\&\quad=(-1)^n
W_n\left(\left(ix+i\frac{\gamma+\delta+1}{2}\right)^2;
\frac{\gamma+\delta+1}{2},\alpha-\frac{\gamma+\delta-1}{2},
\beta+\frac{-\gamma+\delta+1}{2},\frac{\gamma-\delta+1}{2}
\right),
\end{array}
$$
where $\alpha+1=-N$ or $\beta+\delta+1=-N$ or $\gamma+1=-N$,
with $N$ a non-negative integer, and
$$
W_n(x^2;a,b,c,d)\!=\!(-1)^nR_n\left(\lambda(-a\!-\!ix);a+
b\!-\!1,c+d\!-\!1,a+d\!-\!1,a\!-\!d\right),
$$
where $\lambda(t)=t(t+2a)$.

On the other hand taking, for instance,  $\alpha+1=-N$
we get the following factorization for $n>N$:
$$
\begin{array}{rl}
\displaystyle R_n(\lambda(x);\alpha,\beta,\gamma,\delta)=&
\displaystyle R_{N+1}(\lambda(x);-N-1,\beta,\gamma,\delta)
(-1)^{n-N-1}\\ &\hspace{-34mm}\times W_{n-N-1}\left(
\left(i\bigg(x+\frac{\gamma+\delta+1}{2}\bigg)\right)^2;N+
\frac{\gamma+\delta+3}{2},\frac{-\gamma-\delta+1}{2},
\beta+\frac{-\gamma+\delta+1}{2},\frac{\gamma-\delta+1}{2}
\right).
\end{array}
$$

In this case, the Racah polynomials satisfy the following
$\Delta$-Sobolev orthogonality:
$$
\langle p,q\rangle_S=\langle p,q\rangle_d+\left\langle
\left(\frac{\Delta}{\Delta \lambda}\right)^{N+1}p,
\left(\frac{\Delta}{\Delta \lambda}\right)^{N+1}q\right
\rangle_c,
$$
with
$$\begin{array}{rl}
\langle p,q\rangle_d&=\displaystyle\sum_{x=0}^{N}p(x)q(x)
\frac{(\alpha+1,\beta+\delta+1,\gamma+1,\gamma+\delta+1,
(\gamma+\delta+3)/2)_x} {(-\alpha+\gamma+\delta+1,-\beta+
\gamma+1,(\gamma+\delta+1)/2,\delta+1,1)_x},\\
\langle p,q\rangle_c&=\displaystyle\int_{C} p(z^2)q(z^2)
\nu(z i+i+i\racion{\gamma+\delta+N}{2})\nu(-(z
i+i+i\racion{\gamma+\delta+N}{2}))dz,
\end{array}$$
where
$$
\nu(z)\equiv \nu(z;a,b,c,d)=\frac{\Gamma(a+iz)\Gamma(b+iz)
\Gamma(c+iz)\Gamma(d+iz)}{\Gamma(2z)},
$$
being
{\small
$$
a=1+\frac{\gamma+\delta+N}{2},\ b=\frac{-\gamma-\delta-N}
{2}, c=\beta+1+\frac{-\gamma+\delta+N}{2},\ d=1+
\frac{\gamma-\delta+N}{2},
$$
}
and $C$ is the imaginary axis deformed so as
to separate the increasing sequences of poles
$$
k, \ \ -1-\gamma-\delta-N+
k, \ \ \beta-\gamma+k,
\ \ -\delta+k, \qquad k=0, 1, 2, \dots
$$
from the decreasing sequences
$$
-\gamma-\delta-N-2-k, \ \
-1-k, \ \ -\beta-\delta-N
-2-k, \ \ -\gamma-N-2-k,
\qquad k=0, 1, 2, \dots
$$
Of course, we need to assume that these two sets of
poles are disjoint, i.e.,
$$
2a,\ a+b, \ a+c, \ \dots ,\ c+d,\ 2d \not \in
\{0, -1, -2, \dots\}.
$$

On the other hand, in this case, i.e. $\alpha+1=-N$,
we get the following generating functions
(see \cite[p. 29]{kost}) which are valid for all
$x\in \mathbb{C}$:
{\small $$
\begin{array}{rl}
\displaystyle
\sum_{n=0}^\infty \frac{(\alpha+1,\gamma+1)_n}{(\alpha-
\delta+1)_n n!}R_n(\lambda(x);\alpha,\beta,\gamma,\delta)
t^n=& {}_2F_1\left.\left(\begin{array}{c}-x, -x+\beta-
\gamma\\ \beta+\delta+1\end{array}\right| t\right)\\ &
\displaystyle \times {}_2F_1 \left.\left(\begin{array}{c}
x+\alpha+1,x+\gamma+ 1\\ \alpha-\delta+1 \end{array}\!\!
\right| t\right), \\ \displaystyle \sum_{n=0}^\infty
\frac{(\alpha+1,\beta+\delta+1)_n}{(\alpha-\beta-\gamma+
1)_n n!}R_n(\lambda(x);\alpha,\beta,\gamma,\delta)t^n= &
\displaystyle {}_2F_1\hspace{-2mm}\left.\left(
\begin{array}{c}-x, -x-\delta\\ \gamma+1\end{array}\!\!
\right| t\right)\\ & \hspace{-10mm}\displaystyle
\times {}_2F_1 \left.\left(\begin{array}{c}x+\alpha+1,x+
\beta+\delta+ 1\\ \alpha+\beta-\gamma+1 \end{array}\!\!
\right| t\right).
\end{array}
$$
}
See e.g. \cite{wil1} or \cite{kamc} to get more information
about algebraic properties and applications for Wilson
polynomials.
\section{Dual Hahn polynomials}

We can apply an analogous process for the dual
Hahn polynomials
\begin{equation}\label{ex-Hahn}
R_n(\lambda(x);\gamma,
\delta,N)=(\gamma+ 1,-N)_n\,{}_3F_2\left.\left(
\begin{array}{c}-n,-x, x+\gamma+ \delta+1\\\gamma+1,-N
\end{array}\right|1\right),
\end{equation}
with
$\lambda(x)=x(x+\gamma+\delta+1)$, by using the continuous
dual Hahn polynomials \cite[p. 31]{kost}
$$
S_n(x^2;a,b,c)=(-1)^n(a+b,a+c)_n\, {}_3F_2\left.
\left(\begin{array}{c}-n,a+ix,a-ix\\ a+b,a+c\end{array}
\right|1\right).
$$
In fact,
$$ \begin{array}{rl}
&R_n(\lambda(x);\gamma,\delta,N)\\&\quad=(-1)^nS_n
\left(\left(ix+i\frac{\gamma+\delta+1}{2}\right)^2;
\frac{\gamma+\delta+1}{2},\frac{\gamma-\delta+1}{2},
-N-\frac{\gamma+\delta+1}{2}\right),
\end{array}$$
and
$$
S_n(x^2;a,b,c)=(-1)^nR_n\left(\lambda(-a-ix);a+b-1,a-b,
-a-c\right),
$$
where $\lambda(t)=t(t+2a)$.

We get the following factorization for $n>N$:
$$\begin{array}{rl}
\displaystyle R_n(\lambda(x);\gamma,\delta,N)=&
\displaystyle R_{N+1}(\lambda(x); \gamma,\delta,N)~(-1)^{n
-N-1}\\& \times S_{n-N-1}\left(\left(xi+i\frac{\gamma+
\delta+1}{2}\right)^2; N+\frac{\gamma+ \delta+2}{2},\frac{-
\gamma-\delta+1}{2},\frac{\gamma-\delta+1}{2}\right).
\end{array}$$

Dual Hahn polynomials polynomials satisfy the following
$\Delta$-Sobolev orthogonality:
$$
\langle p,q\rangle_S=\langle p,q\rangle_d+\left\langle
\left(\frac{\Delta}{\Delta \lambda}\right)^{N+1}p,\left(
\frac{\Delta}{\Delta\lambda}\right)^{N+1}q\right\rangle_c,
$$
with
$$\begin{array}{rl}
\langle p,q\rangle_d&=\displaystyle \sum_{x=0}^{N}p(x)q(x)
\frac{(2x+\gamma+\delta+1)(\gamma+1,-N)_x}{(-1)^x(x+\gamma+
\delta+1)_{N+1}(\delta+1,1)_x},\\ \langle p,q\rangle_c&
=\displaystyle \int_{C} p(z^2)q(z^2)w(z;\gamma,\delta,N)dz,
\end{array}$$
where
$$
w(z;\alpha,\beta,N)=\nu(z
i+i+i\racion{\gamma+\delta+N}{2};\gamma,\delta,N)\nu(-zi
-i-i\racion{\gamma+\delta+N}{2};\gamma,\delta,N),
$$
being
$$
\nu(z;a,b,c)=\frac{\Gamma(a+iz)\Gamma(b+iz)\Gamma(c+iz)}
{\Gamma(2zi)},
$$
and
$$
a=1+\frac{\gamma+\delta+N}{2}\,,\qquad
b=1+\frac{\gamma-\delta+N}{2}\,,\qquad
c=-\frac{\gamma+\delta+N}{2},
$$
and $C$ is the imaginary axis deformed so as to separate
the increasing sequences of poles
$$
k, \quad -\gamma-
\delta-N-1+k,  \quad
-\delta+k, \qquad k=0, 1, 2, \dots,
$$
from the decreasing sequences
$$
-\gamma-\delta-N-2-k,
\quad -1-k,\quad -\gamma-
N-2-k, \qquad k=0, 1, 2, \dots
$$
Of course, we need to assume that these
two sets of poles are disjoint, i.e.,
$$
2a,\ a+b, \ a+c, \ \dots ,\ 2c\not \in \{0, -1, -2, \dots\}.
$$

On the other hand we get the following generating functions
(see \cite[p. 36]{kost}) which are valid for all $x\in
\mathbb{C}$:
$$
(1-t)^{N-x}{}_2F_1\left.\left(\begin{array}{c} -x,-x-\delta
\\ \gamma+1 \end{array}\right|t\right)=\sum_{n=0}^\infty
\frac{(-N)_n}{n!}R_n(\lambda(x); \gamma,\delta,N)t^n.
$$
$$
(1-t)^x {}_2F_1\left.\left(\begin{array}{c} x-N,x+\gamma+1\\
-\delta-N\end{array}\right|t\right)=\sum_{n=0}^\infty\frac{(
\gamma+1,-N)_n}{(-\delta-N)_n n!}R_n(\lambda(x);\gamma,
\delta,N)t^n.
$$

\section{Krawtchouk polynomials}\label{SecKrawtchouk}

Similarly, properties for the Krawtchouk polynomials
\begin{equation}
\label{ex-Krawtchouk}
K_n(x;p,N)=(-N)_n p^n\,{}_2F_1\left(\left. \begin{array}{cc}
-n, -x \\-N\end{array} \right|\frac 1p \right),
\end{equation}
with $p\in \mathbb{C}$, $p\ne 0, 1$, can be obtained via
Meixner polynomials
$$
M_n(x;\beta,c)=\frac{c^n(\beta)_n}{(c-1)^n}\, {}_2F_1\left(
\left.\begin{array}{c}-n, -x \\ \beta \end{array}\right|1-
\frac 1c\right),\quad \beta>0, \ 0<c<1.
$$
In fact,
$$\begin{array}{rl}
K_n(x;p,N)&=M_n\left(x;-N,\frac{p}{p-1}\right),\\
M_n(x;\beta,c)&=K_n\left(x;-\beta,\frac{c}{c-1}\right),
\end{array}$$
setting $\beta =-N$ and $c=p/(p-1)$.

We have the following factorization for the Krawtchouk
polynomials for $n>N$:
$$
\begin{array}{rl}
K_n(x;p,N)=&  \displaystyle K_{N+1}(x;p,N) M_{n-N-1}
(x-N-1;N+2,p/(p-1)).
\end{array}
$$
Furthermore, these polynomials satisfy the following
$\Delta$-Sobolev orthogonality:
$$
\langle r,q\rangle_S=\langle r,q\rangle_d+\left\langle
\Delta^{N+1}r,\Delta^{N+1}q\right\rangle_c,
$$
with
\begin{eqnarray}
\langle r,q\rangle_d&\displaystyle =\sum_{x=0}^{N}r(x)q(x)
\frac{p^x (1-p)^{N-x}}{\Gamma(x+1)\Gamma(N-x+1)},
\label{kraw-ort}\\
\langle r,q\rangle_c&=\displaystyle \int_C r(z)q(z)
\Gamma(-z)\Gamma(1+z)\left(\frac{p}{1-p}\right)^z dz,
\end{eqnarray}
where  $C$ is the imaginary axis deformed so as to separate
the increasing from the decreasing sequence of poles of the
weight function, in fact in this case we can consider the
curve $C=\{-1/2+ti:t\in\mathbb{R}\}$.

\begin{rem} Note that the property of orthogonality for the
Krawtchouk polynomials (\ref{kraw-ort}) is valid for all
$p\in\mathbb{C}$, with $p\ne 0, 1$ by using an analytic
continuation for the standard weight function associated
with the Krawtchouk polynomials.
\end{rem}

On the other hand, we get the following generating function
(see \cite[p. 47]{kost}) which is valid for all $x\in
\mathbb{C}$:
$$
\left(1-\frac{1-p}{p}~t\right)^x(1+t)^{N-x}= \sum_{n=
0}^\infty {N \choose n} K_n(x;p,N)t^n.
$$
\section{Limit relations between hypergeometric orthogonal
polynomials}
In this section, we study the limit relations involving
the orthogonal polynomials, considered in this paper,
associated with some families of polynomials of the
Askey scheme of hypergeometric orthogonal polynomials
\cite{kost}.

Let us now consider such limits for any $n\in \mathbb{N}_0$:
\begin{enumerate}
\item Racah $\to$ Hahn. If we take $\gamma+1=-N$ and
$\delta \to \infty$ in the definition (\ref{ex-Racah}) of
the Racah polynomials, we obtain the Hahn polynomials
defined by (\ref{HahnHipergeometric}).
Hence
$$
\lim_{\delta\to \infty} R_n(\lambda(x);\alpha,\beta,-N-1,
\delta)= h_n^{\alpha,\beta}(x;N).
$$
The Hahn polynomials can also be obtained from the Racah
polynomials by taking $\delta=-\beta-N-1$ in the definition
(\ref{ex-Racah}) and letting $\gamma\to \infty$:
$$
\lim_{\gamma\to \infty} R_n(\lambda(x);\alpha,\beta,\gamma,
-\beta-N-1))=h_n^{\alpha,\beta}(x;N).
$$
Another way to do this is to take $\alpha+1=-N$ and $\beta
\to \beta+\gamma+N+1$ in the definition (\ref{ex-Racah}) of
the Racah polynomials and then take the limit $\delta\to
\infty$.
In that case we obtain the Hahn polynomials
given by (\ref{HahnHipergeometric}) in the following way:
$$
\lim_{\delta\to\infty} R_n(\lambda(x);-N-1,\beta+\gamma+N+1,
\gamma,\delta)=h_n^{\gamma,\beta}(x;N).
$$
\item  Racah $\to$ Dual Hahn. If we take $\alpha+1=-N$ and
let $\beta \to \infty$ in the definition (\ref{ex-Racah})
of the Racah polynomials, then we obtain the dual Hahn
polynomials defined by (\ref{ex-Hahn}).
Hence
$$
\lim_{\beta\to \infty} R_n(\lambda(x);-N-1,\beta,\gamma,
\delta)= R_n(\lambda(x);\gamma,\delta,N).
$$
If we take $\beta=-\delta-N-1$ and $\alpha\to \infty$ in
(\ref{ex-Racah}), then we also obtain the dual Hahn
polynomials:
$$
\lim_{\alpha\to \infty} R_n(\lambda(x);\alpha,-\delta-N-1,
\gamma, \delta)=R_n(\lambda(x);\gamma,\delta,N).
$$
Finally, if we take $\gamma+1=-N$ and $\delta\to \alpha +
\delta+N+1$ in the definition (\ref{ex-Racah}) of the Racah
polynomials and take the limit $\beta\to\infty$ we find the
dual Hahn polynomials given by (\ref{ex-Hahn}) in the
following way:
$$
\lim_{\beta\to \infty} R_n(\lambda(x);\alpha,\beta,
-N-1,\alpha +\delta+N+1)=R_n(\lambda(x);\alpha,\delta,N).
$$
\item Hahn $\to$ Krawtchouk.
If we take $\alpha = (1-p)t$ and $\beta = pt$ in the
definition (\ref{HahnHipergeometric}) of the Hahn
polynomials and let $t\to \infty$ we obtain the Krawtchouk
polynomials defined by (\ref{ex-Krawtchouk}):
$$
\lim_{t\to \infty} h_n^{(1-p)t,pt}(x;N)= K_n(x; p,N).
$$
\item Dual Hahn $\to$ Krawtchouk.
In the same way we find the Krawtchouk polynomials from the
dual Hahn polynomials by setting $\gamma= pt$,
$\delta=(1-p)t$ in (\ref{ex-Hahn}) and letting $t\to \infty$:
$$
\lim_{t\to \infty} R_n(\lambda(x);pt,(1-p)t,N)=K_n(x; p,N).
$$
\end{enumerate}

\begin{rem}
The proof of each one of these limits is straightforward
once one reduces the hypergeometric representation of
each family as we did for the Hahn polynomials (see
(\ref{HahnHypergeometricreduced})) which is valid for
all $n\in \mathbb{N}_0$.
\end{rem}

\appendix
\section{Orthogonality relations for Meixner polynomials
with general parameter}

In this appendix we will show that Meixner polynomials,
$\{M_n(x;\beta,c)\}_{n\ge 0}$, with  $c<0$ and
$\beta\in\mathbb{C}$, can be provided with a property of
orthogonality which can be obtained through a process
limit from the continuous Hahn polynomials.

From the TTRR of the continuous Hahn polynomials it is
straightforward to obtain the following relation
\begin{equation} \hspace{-7mm}
p_{n+1}(ix;0, -\racion tc, t, \beta)=(x\!-\!B_n)
 p_n(ix;0, -\racion tc, t, \beta)\!-\!C_n  p_{n-1}(ix;0,
 -\racion tc, t, \beta),\label{TTRR_CH_to_Meixner}
\end{equation}
where $t\in \mathbb{R}^+$, and
$$
\begin{array}{rl}
B_n&= \frac{n\left(n-\frac{t}{c}+t-1\right)\left(n-\frac{t}
{c}+\beta-1\right)}{\left(2n-\frac{t}{c}+t+\beta-2\right)
\left(2n-\frac{t}{c}+t+\beta-1\right)}- \frac{\left(n-
\frac{t}{c}+t+\beta-1\right)(n+t)(n+\beta)}{\left(2n-\frac
{t}{c}+t+\beta-1\right)\left(2n-\frac{t}{c}+t+\beta
\right)},\\ C_n&=\frac{n \left(n-\frac{t}{c}+t-1\right)
\left(n-\frac{t}{c}+\beta-1\right)\left(n-\frac{t}{c}+t+
\beta-2\right)(n+t-1)(n+\beta-1)}{\left(2n-\frac{t}{c}+t+
\beta-1\right)\left(2n-\frac{t}{c}+t+\beta-2\right)^2\left(
2n-\frac{t}{c}+t+\beta-3\right)}.
\end{array}
$$
Since
$$\begin{array}{rl}
\lim_{|t|\to\infty}B_n&=\displaystyle \frac{n+c(n+\beta)}
{1-c},\\ \lim_{|t|\to\infty}-C_n&=\displaystyle \frac{nc(n+
\beta-1)}{(c-1)^2}, \end{array}
$$
coincide with the coefficients of the TTRR for the monic
Meixner polynomials, with initial conditions $p_{-1}=0$ and
$p_0=1$, one deduces
\begin{equation}\label{CHtoMeixner}
\lim_{|t|\to\infty}(-i)^n p_n(ix;0, -t/c, t, \beta)=
M_n(x;\beta,c),\qquad n=0, 1, 2, \dots,
\end{equation}
by using induction in (\ref{TTRR_CH_to_Meixner}).

\begin{prop} For any $\beta, c\in\mathbb{C}$,
$c\notin[0,\infty)$ and $-\beta\not \in \mathbb{N}$,
the following property of orthogonality for the Meixner
polynomials fulfills:
\begin{equation}
\label{Ortog Meixner}\hspace{-7mm}\int_{C}\!\!M_n(z;c,\beta)
z^m\Gamma(-z)\Gamma(\beta+z)(-c)^z dz=0,\ 0\le m<n, \
n=0, 1, 2, \dots
\end{equation}
where $C$ is a complex contour from $-\infty i$ to
$\infty i$ separating  the increasing poles $\{0, 1, 2,
\dots \}$ from the decreasing poles $\{-\beta, -\beta-1,
-\beta-2, \dots\}$.
\end{prop}
\begin{pf}
We prove the result for $c<0$, thus the general case is
obtained by analytic continuation.

Let us assume that $\beta$ is such that the contour
$C=\{-1/2+y i: y\in\mathbb{R}\}$ separates the poles of
$\Gamma(\beta+z)$ from the poles of $\Gamma(-z)$, i.e.,
$\Re \beta> 1/2$, and let us take the normalized weight
for the continuous Hahn polynomials (see (\ref{weigthch}))
$$
W_t(z)=\frac{w(iz;0,-t/c,t,\beta)}{\Gamma(-t/c)\Gamma(t)}=
\Gamma(-z)~\frac{\Gamma(-t/c-z)}{\Gamma(-t/c)}
~\frac{\Gamma(t+z)}{\Gamma(t)}~\Gamma(\beta+z).
$$
Notice that
$$
\lim_{t\to \infty} W_t(z)=\Gamma(-z)
\Gamma(\beta+z)(-c)^z=:W(z),
$$
pointwise in $C$ by using the Stirling formula
$$
\Gamma(z)=\sqrt{2\pi}z^{z-\half}e^{-z}(1+o(1)),\qquad z\to
\infty,|\arg(z)|<\pi.
$$

It is known that
$$
|\Gamma(x+i y)|\leq|\Gamma(x)|,\qquad \forall\, x,y\in
\mathbb{R},
$$
hence
\begin{equation}\label{Dominatedweight}
|W_t(z)|\leq \left|\Gamma(-z)~\Gamma(\beta+z)\right|\,.
\end{equation}
 Using once again induction on
(\ref{TTRR_CH_to_Meixner}) and due to the exponential
behavior of the right-hand side of (\ref{Dominatedweight})
at the endpoints of $C$, one obtains that
$$
p_n(iz;0,-t/c,t,\beta)W_t(z)
$$
is dominated
by an integrable function on $C$. Thus, from the dominated
convergence theorem
$$
\lim_{t\to \infty}\int_C (-i)^{n+1}p_{n}(iz;0, -t/c, t,
\beta)z^m W_t(z)dz=\int_C M_n(z;\beta,c) z^m W(z) dz.
$$
On the other side since $C$ also separates the poles of
$\Gamma(-t/c-z)$ from the poles of $\Gamma(t+z)$, we get
$$
\int_C (-i)^{n+1}p_{n+1}(iz;0, -t/c, t,\beta)z^m W_t(z)=0.
$$
Thus (\ref{Ortog Meixner}) holds for $\Re \beta>1/2$.

The general case is straightforward by using that if $C$
is a contour separating the poles and $C_1=\{\lambda+yi:
y\in \mathbb{R}\}$, where $\lambda\in (0,1)$, $\lambda\ne
\Re \beta$, which does not separate the poles, then the
integral through $C$ and $C_1$ differs on a finite number
of residues.
\qed \end{pf}

The case $c>0$ cannot be considered by an integral of the
form (\ref{Ortog Meixner}) since it diverges. However,
when $|c|<1$, (\ref{Ortog Meixner}) is rewritten on the
form (see \cite[\S 5.6]{Nico} for details)
$$
\sum_{x=0}^\infty M_n(x;c,\beta) x^m~\frac{\Gamma(\beta+
x)~c^x}{x!}=0,
$$
which is also valid for $c\in(0,1)$ and coincides with the
very well-known orthogonal relations for Meixner
polynomials.

{\bf Acknowledgements:}
We thank referees for their suggestions which have improved
the presentation of the paper.
\\
The authors also wish to thank  R. \'{A}lvarez-Nodarse,
F. Marcell\'{a}n, J.J. Moreno-Balc\'{a}zar and A. Zarzo, for their
useful suggestions and comments.

\end{document}